\documentclass[12pt]{article}
\newtheorem{prop}{Proposition}[section]
\newtheorem{lem}[prop]{Lemma}

\newtheorem{thm}[prop]{Theorem}
\usepackage{latexsym}
\usepackage{amssymb}
\usepackage[dvips]{color}
\usepackage{graphics}
\usepackage{graphicx}
\usepackage{pstricks}

\newcommand{\covb}{\prec}
\newcommand{\covr}{\succ}

\begin{document}

\bibliographystyle{elsart-num-sort}
\title{The M\"{o}bius Function\\ of a\\ Restricted Composition Poset}

\author{Adam M. Goyt\footnote{This work was partially done while the author was visiting DIMACS.}\\ Department of Mathematics\\ Minnesota State University Moorhead\\ 1104 7th Avenue South\\ Moorhead, Minnesota 56563\\ goytadam@mnstate.edu\\www.mnstate.edu/goytadam}

\maketitle


\medskip


\begin{abstract}

We study a poset of compositions restricted by part size under a partial ordering introduced by Bj\"{o}rner and Stanley.  We show that our composition poset $C_{d+1}$ is isomorphic to the poset of words $A_d^*$.  This allows us to use techniques developed by Bj\"{o}rner to study the M\"{o}bius function of $C_{d+1}$.  We use counting arguments and shellability as avenues for proving that the M\"{o}bius function is $\mu(u,w)=(-1)^{|u|+|w|}{w\choose u}_{dn}$, where ${w\choose u}_{dn}$ is the number of $d$-normal embeddings of $u$ in $w$.  We then prove that the formal power series whose coefficients are given by the zeta and the M\"{o}bius functions are both rational.  Following in the footsteps of Bj\"{o}rner and Reutenauer and Bj\"{o}rner and Sagan, we rely on definitions to prove rationality in one case, and in another case we use finite-state automata.  

\end{abstract}

\section{Introduction}

Let $\mathbb{P}=\{1,2,3,\dots\}$, $[n]=\{1,2,3,\dots,n\}$ and $[i,j]=\{i,i+1,i+2,\dots,j\}$.  A {\it composition} $\alpha$ is an element $(\alpha_1,\alpha_2,\dots,\alpha_k)\in \mathbb{P}^k$.  The {\it length} of a composition $\alpha=(\alpha_1,\alpha_2,\dots,\alpha_k)$ is $\ell(\alpha)=k$, and the {\it norm}, $|\alpha|=n$, is $n$ if $\sum_{i=1}^k\alpha_i=n$.  The $\alpha_i$ are called {\it parts}.  If $|\alpha|=n$, then we say that $\alpha$ is a composition of $n$.  Define $C_n$ to be the set of all compositions of $n$ and let $$C=\bigcup_n C_n.$$ 

The partial ordering on $C$ that we are interested in was introduced by Bj\"{o}rner and Stanley and can be found in~\cite{RStanleyVol1Ed2}. To define the poset $C$ we need only define its covering relation $\alpha \covb \beta$.  In $C$ we say that $\beta\covr\alpha=(\alpha_1,\alpha_2,\dots,\alpha_k)$ if $\beta$ is of one of the two following forms:
\begin{eqnarray*}\beta&=&(\alpha_1,\alpha_2,\dots,\alpha_{i-1},\alpha_i+1,\alpha_{i+1},\dots,\alpha_k)\\ 
\beta&=&(\alpha_1,\alpha_2,\dots,\alpha_{i-1},\alpha_i+1-h,h,\alpha_{i+1},\dots,\alpha_k),\\
\end{eqnarray*} where $h\leq \alpha_i$.

We will consider the subposet $C_d$, which is the set of compositions with part sizes at most $d$, that is $\alpha=(\alpha_1,\alpha_2,\dots,\alpha_k)$ is an element of $C_d$ if and only if $\alpha_i\leq d$ for $1\leq i\leq k$.  We will show that this poset is isomorphic to a poset of words in order to determine the M\"{o}bius function of $C_d$. 

Let $A^*$ be the free monoid under concatenation of $A=\{a,b\}$.  We think of $A^*$ as the set of all words that can be created from the {\it alphabet} $A$.  The identity in $A^*$ is the empty word $\epsilon$.    We say that the length of a word $u=u_1u_2\dots u_k$ is $|u|=k$.  Let  $u=u_1u_2\dots u_k$ and $w=w_1w_2\dots w_\ell$ be words in $A^*$.  We make $A^*$ into a poset by letting $u\leq w$ if there exist $i_1\leq i_2\leq\dots\leq i_k$ such that $u_j=w_{i_j}$ for $1\leq j\leq k$.  We call the set $\iota=\{i_1,i_2\dots,i_k\}$ an {\it embedding} of $u$ in $w$ and let $w_\iota=w_{i_1}w_{i_2}\dots w_{i_k}$.  If $\iota$ is an embedding of $u$ in $w$, then we say that $w_j$ is {\it supported} by $u$ in $\iota$ if $j\in\iota$.  For example, the word $abaab$ is a subword of $w=aabbababb$, since $w_2w_3w_5w_7w_8=abaab$, and $w_3$ is supported in this embedding.  

Now, let $\phi:\mathbb{P}\rightarrow A^*$ be given by $\phi(k)=a\underbrace{bb\dots b}_{k-1}$.  Given any composition $\alpha=(\alpha_1,\alpha_2,\dots,\alpha_k)\in C$, , let $\phi(\alpha)=\phi(\alpha_1)\phi(\alpha_2)\cdots\phi(\alpha_k)$ omitting the initial $a$ from $\phi(\alpha_1)$.  Clearly, $\phi:C\rightarrow A^*$ is a bijection.  
   
\begin{lem} The map $\phi$ is an isomorphism of $A^*$ and $C$ as partially ordered sets. \end{lem}

\noindent {\bf Proof:} Since we have already established $\phi$ as a bijection between $A^*$ and $C$, it is enough to show that $\phi$ preserves the partial orderings.  

Let $\alpha=(\alpha_1,\alpha_2,\dots,\alpha_k)\covb\beta=(\beta_1,\beta_2,\dots,\beta_\ell)\in C$.  If $\ell=k$ then, for some $j$, $\beta_j=\alpha_j+1$ and $\beta_i=\alpha_i$ for $i\not=j$.  In this case $\phi(\beta)$ is obtained from $\phi(\alpha)$ by inserting a $b$ into $\phi(\alpha)$ anywhere between the $(j-1)^{st}$ and $j^{th}$ occurrence of an $a$.  Thus, $\phi(\alpha)$ is a subword of $\phi(\beta)$.

If $\ell=k+1$ then $\beta_j=\alpha_j+1-h$ for some $h\leq \alpha_j$, $\beta_{j+1}=h$ and $\beta_i=\alpha_i$ for $i\not=j$ or $j+1$.  In this case $\phi(\beta)$ is obtained from $\phi(\alpha)$ by inserting an $a$ between the $(j-1)^{st}$ and $j^{th}$ occurrence of an $a$, so that $h-1$ $b$'s follow the inserted $a$.

Thus $\beta\covr\alpha$ implies $\phi(\beta)\covr\phi(\alpha)$.  The converse is proved similarly.  $\square$

\medskip

Let $A^*_d$ be the subposet of $A^*$ consisting of all words that do not have $d+1$ consecutive $b$'s.  Notice that $\phi$ restricts to an isomorphism between $C_{d+1}$, as defined above, and the subposet, $A^*_d$.  Much is known about partially ordered sets on words and subword order, so it will be convenient to work with the poset $A^*_d$ to understand the poset $C_{d+1}$.

We assume familiarity with the M\"{o}bius function of a locally finite poset~\cite{RStanleyVol1}.  To understand the M\"{o}bius function of the poset $A_d^*$ we will adapt the M\"{o}bius function for $A^*$.  In~\cite{Bjornersubword}, Bj\"{o}rner shows that the M\"{o}bius function for $A^*$ with the subword order described above is given by $$\mu(u,w)=(-1)^{|u|+|w|}{w\choose u}_n,$$ where ${w\choose u}_n$ is the number of normal embeddings of $u$ in $w$.  We will alter his definition of normal embeddings to define $d$-normal embeddings and show that in the poset $A^*_d$ we have the following.

\begin{thm} For elements $u,$ $w\in A^*_d$, we have $$\mu(u,w)=(-1)^{|u|+|w|}{w\choose u}_{dn},$$ where ${w\choose u}_{dn}$ is the number of $d$-normal embeddings of $u$ in $w$.  \end{thm}

We now define a $d$-normal embedding and right-most embedding.  We will call a set $[r,s]$ a {\it run} if $w_r=w_{r+1}=\dots=w_s$.  Let the {\it repetition set} of $w$ be $R(w)=\{j:w_j=w_{j-1}\}$.  An embedding $\iota=\{i_1,i_2,\dots,i_k\}$ of $u$ in $w$ is called {\it $d$-normal} if 

\begin{enumerate}
\renewcommand{\labelenumi}{\alph{enumi}}
\renewcommand{\labelenumii}{\roman{enumii}}
\item $R(w)\subseteq\iota$, and 
\item if $u$ has a run of $d$ $b$'s and $w_{i_j}$ corresponds to the first $b$ in this run then $i_j=1$ or $w_{i_j-1}=a$ and $i_j-1\in\iota$ or $i_j-1=1$.
\end{enumerate}

This differs from a normal embedding in that a normal embedding must only satisfy part a of the definition of a $d$-normal embedding.  

The {\it right-most embedding} of $u$ in $w$ is the unique embedding $\{j_1,j_2\dots j_r\}$ such that for any other embedding $\{i_1,i_2\dots i_r\}$ of $u$ in $w$ we have that $i_\ell\leq j_\ell$, for each $1\leq \ell\leq r$.  

The definition of a $d$-normal embedding, while somewhat cumbersome, is as simple as possible and is borne out of the proof involving shellability of the poset.  For an explanation of how the definition was borne, please see the end of Section 3.

An example is in order to help describe the definition of $d$-normal.  Consider the poset $A^*_3$ and the interval of this poset $[bbb,bbabb]$.  Let $u=bbb$ and $w=bbabb$.  In any $d$-normal embedding of $u$ in $w$ the $b$'s located in positions 2 and 5 must be supported according to condition a.  This means that the embeddings $\{1,2,5\}$ and $\{2,4,5\}$ are the only possible $d$-normal embeddings.  Now, $u$ contains a run of three $b$'s, so by a condition b $\{2,4,5\}$ cannot be a 3-normal embedding because $w_2$ corresponds to the first $b$ in $u$, and $w_1\not=a$.  Thus, $\{1,2,5\}$ is the only $d$-normal embedding of $u$ in $w$.  So $\mu(u,w)=(-1)^{3+5}\cdot1=1$.  A quick sketch of the interval $[bbb,bbabb]$ in $A^*_3$ confirms this.

In Section 2 we prove Theorem 1.2 using bijective techniques, and in Section 3 we prove it using shellability.  The proofs in Sections 2 and 3 are adaptations of similar proof given by Bj\"{o}rner.  In Section 4 we show that the zeta function and M\"{o}bius function are rational.  Finally, in the last section we determine a generating function in commuting variables for the zeta function and the M\"{o}bius function.

\section{Combinatorial Proof}

We remind the reader that the M\"{o}bius function of a poset is the unique function satisfying the following.
\begin{itemize}
\item $\mu(x,x)=1$
\item $\sum_{x\leq z\leq y}\mu(z,y)=0$
\item $\mu(x,y)=0$ for $x\nleq y$
\end{itemize}

\noindent {\bf Proof of Theorem 1.2:} This proof is essentially the same as the proof of Theorem 3.1 from~\cite{Bjornerfactor}.  We will show that the function $f(u,w)=(-1)^{|w|+|u|}{w\choose u}_{dn}$ satisfies the three conditions above and hence, must be $\mu$.  

Clearly, $f(u,u)=1$, and $f(u,w)=0$ for $u\nleq w$.  It remains to show that\\ $\sum_{u\leq v\leq w}f(v,w)=0$ for any $v\leq w$.  

Suppose that $|w|=n$ and $|u|=k$.  Let $N$ be the set of embeddings $\iota$ such that $\iota$ is a $d$-normal embedding of $w_\iota$ in $w$, and $u\leq w_\iota$.  Let $N_e=\{\iota\in N:|\iota|\mathrm{\:is\:even}\}$ and $N_o=\{\iota\in N:|\iota|\mathrm{\:is\:odd}\}$.  Then $$\sum_{u\leq v\leq w}f(v,w)=(-1)^n(\#N_e-\#N_o).$$  A bijection between $N_e$ and $N_o$ will complete the proof.

 For each embedding $\iota\in N$ let $\iota_f$ be the minimal element of $[n]$ such that $\iota_f$ is not in the right-most embedding of $u$ in $w_\iota$.  Define $\psi:N\rightarrow N$ by $$\psi(\iota)=\left\{\begin{array}{ll} \iota\cup\{\iota_f\}&\mathrm{if\:}\iota_f\notin\iota,\\ \iota-\{\iota_f\}&\mathrm{if\:}\iota_f\in\iota.\end{array}\right.$$  

We will show that $\psi$ is well-defined.  If $\iota_f\notin\iota$ then, since $R(w)\subseteq\iota$, we must have $R(w)\subseteq\iota\cup\{\iota_f\}$.  If $\iota_f\in\iota$ then we want to show that $\iota_f\notin R(w)$.  If $\iota_f\in R(w)$ then $w_{\iota_f}=w_{\iota_f-1}$.  Also $\iota_f-1$ must be in the right-most embedding of $u$ in $w_\iota$.  This is impossible by the definition of right-most embedding.  

We must show that $w_{\psi(\iota)}$ is still an element of $A^*_d$.  If $\psi$ removes an $a$ creating a run of $b$'s of length $d+1$  in $w_{\psi(\iota)}$ then the first $b$ in the run of $b$'s immediately preceding $a$ cannot be in the right-most embedding of $u$ in $w_{\iota}$.  Thus, the $a$ would never have been removed by $\psi$.  

If $\psi$ inserts a $b$ creating a run of $b$'s of length $d+1$ in $w_{\psi(\iota)}$ then $w_{\iota}$ has a run of $b$'s of length $d$.  Let the first $b$ of this run be in position $i_t$.  If $i_t=1$ then a $b$ must have been inserted at the beginning of $w_\iota$, which is impossible by part b of the definition of $d$-normal.  If $i_t\not=1$, then, by part b of the definition of $d$-normal, we must have $w_{i_t-1}=a$ and $i_t-1\in{\iota}$ or $\iota_t-1=1$ rendering the creation of this run impossible.

Since $\psi$ is its own inverse and changes the parity of $|\iota|$, the proof is complete. $\square$

\section{Shellability}

We now give an alternate proof that $\mu(u,w)=(-1)^{|w|+|u|}{w \choose u}_{dn}$ by showing that $A^*_d$ is dual CL-shellable.  Throughout this section we assume familiarity with lexicographical shellability.  The same labeling that Bj\"{o}rner used in~\cite{Bjornersubword} to determine the M\"{o}bius function of $A^*$ will work for our poset $A^*_d$.   

Let $[u,w]$ be an interval of $A^*$ and $\mathcal{E}^*(A^*)$ be the set of edges of maximal chains, $m=(u=x_1,x_2,\dots,x_k=w)$ in $[u,w]$.  Let $\ell:\mathcal{E}^*(A^*)\rightarrow \mathbb{Z}$ be defined as follows.  Label the edge $(x_{k-1},x_k,m)$ by $\ell(x_{k-1},x_k,m)=i_1$, where $i_1$ is the smallest element of $[k]$ such that $[k]-\{i_1\}$ is an embedding of $x_{k-1}$ in $x_{k}$.  Label the edge $(x_{k-2},x_{k-1},m)$ by $\ell(x_{k-2},x_{k-1},m)=i_2$, where $i_2$ is the least element of $[k]-\{i_1\}$ such that $[k]-\{i_1,i_2\}$ is an embedding of $x_{k-2}$ in $x_k$.  Repeat this process for the remaining edges.  It is clear that if two maximal chains have the same first $s$ edges then the labels on these edges will be the same.  Figure 1 shows a chain-edge labeling for $[abb,aabbab]$.

\begin{thm}[Bj\"{o}rner] The map $\ell:\mathcal{E}^*(A^*)\rightarrow \mathbb{Z}$ described above is a  dual CL-labeling for each interval of $A^*$.  \end{thm}

In Bj\"{o}rners proof, \cite{Bjornersubword}, he shows that the unique increasing maximal chain is the one corresponding to the right-most embedding of $u$ in $w$.  Let $\{i_1,i_2,\dots,i_k\}$ be the rightmost embedding of $u$ in $w$, then the increasing maximal chain is the maximal chain labeled by the elements of $[n]-\{i_1,i_2,\dots,i_k\}$ from smallest to largest.  The unique increasing chain in Figure 1 is $$abbaba\stackrel{1}{-}abbab\stackrel{3}{-}abab\stackrel{5}{-}abb.$$  It is also clear that since we are working with the right-most embedding, this chain must have the lexicographically smallest possible entries. 

Each maximal chain of an interval $[u,w]$ in $A^*_d$ is also a maximal chain of the same interval in $A^*$.  For $u,w\in A^*_d$, we will denote by $[u,w]_d$ the interval $[u,w]$ in $A^*$ restricted to $A^*_d$.  We will label the edges of $A^*_d$ in exactly the same way as we did for $A^*$.  

\begin{thm} The map $\ell_d:\mathcal{E}^*(A^*_d)\rightarrow \mathbb{Z}$ is a dual CL-labeling for each interval of $A^*_d$, and hence $A^*_d$ is CL-shellable.  \end{thm}
 
\noindent {\bf Proof:}  We already know that there is exactly one increasing maximal chain, $m_0$, in the labeling of an interval $[u,w]$ in $A^*$, and that $m_0$ corresponds to the right-most embedding of $u$ in $w$.  This means that there is at most one increasing maximal chain among those in the interval $[u,w]_d$.  We must show that this chain is indeed in $[u,w]_d$.  Unless a chain from $[u,w]$ passes through an element with a run of $b$'s of length at least $d+1$, that chain is in $[u,w]_d$.  

Suppose that the unique increasing maximal chain for $[u,w]$ passes through an element with a run of $b$'s of length at least $d+1$ and hence is not in $[u,w]_d$.  At some step in the labeling of this chain we must have $x_i\covr x_{i-1}$, where $x_i$ does not have a run of $b$'s of length at least $d+1$, but $x_{i-1}$ does.  This may only happen if an $a$ that was separating those $b$'s in $x_i$ is removed in this step.  Since this chain corresponds to the right-most embedding of $u$ in $w$ and $u$ and $w$ do not have a run of $b$'s of length at least $d+1$, the $a$ removed could not have been the left-most removable element.   $\square$

\vspace{.5in}

\begin{pspicture}(-5,-5)(5,5)
\rput[l](-5,5){\bf Figure 1}
\rput[t](0,4.5){$aabbab$}
\rput[t](-4.5,2){$aabbb$}
\rput[t](-1.5,2){$aabba$}
\rput[t](1.5,2){$abbab$}
\rput[t](4.5,2){$aabab$}
\rput[t](-4.5,-2){$abbb$}
\rput[t](-1.5,-2){$aabb$}
\rput[t](1.5,-2){$abba$}
\rput[t](4.5,-2){$abab$}
\rput[t](0,-4.5){$abb$}
\psline(0,4.2)(-4.5,2)
\psline(0,4.2)(-1.5,2)
\psline(0,4.2)(1.5,2)
\psline(0,4.2)(4.5,2)
\psline(-4.5,1.7)(-4.5,-2)
\psline(-4.5,1.7)(-1.5,-2)
\psline(-1.5,1.7)(-1.5,-2)
\psline(-1.5,1.7)(1.5,-2)
\psline(1.5,1.7)(-4.5,-2)
\psline(1.5,1.7)(1.5,-2)
\psline(1.5,1.7)(4.5,-2)
\psline(4.5,1.7)(-1.5,-2)
\psline(4.5,1.7)(4.5,-2)
\psline(-4.5,-2.3)(0,-4.5)
\psline(-1.5,-2.3)(0,-4.5)
\psline(1.5,-2.3)(0,-4.5)
\psline(4.5,-2.3)(0,-4.5)
\rput[t](-3.2,3){5}
\rput[t](-1.2,3){6}
\rput[t](1.2,3){1}
\rput[t](3.2,3){3}
\rput[t](-4.7,0){1}
\rput[t](-3.3,0){3}
\rput[t](-3.7,-1){5}
\rput[t](-1.7,1){5}
\rput[t](-0.3,0){1}
\rput[t](1.2,0.5){6}
\rput[t](3.7,1){5}
\rput[t](4.8,0){1}
\rput[t](3.1,-0.5){3}
\rput[t](-3.2,-3){3}
\rput[t](-1.2,-3){1}
\rput[t](1.2,-3){5}
\rput[t](3.2,-3){5}
\end{pspicture}

\vspace{.25in}

\noindent {\bf Proof of Theorem 1.2:}  It suffices to show that the number of decreasing chains in $[u,w]_d$ is the same as the number of $d$-normal embeddings of $u$ in $w$.  

Suppose $\iota=\{i_1,i_2\dots, i_k\}$ is an embedding of $u$ in $w$ such that there is a maximal chain $m$ with $[n]-\{\ell_1(m),\ell_2(m),\dots,\ell_{n-k}(m)\}=\iota$, where $\ell_j(m)$ is the label on the $j^{th}$ edge of $m$.  Then $m$ is decreasing if $m$ is obtained by deleting the entries $w_j$, $j\in [n]-\iota$, from right to left.  The only way that this is possible is if for each $j\in[n]-\iota$, we have that $w_j\not=w_{j-1}$.  This means that the embedding satisfies part a of the definition of a $d$-normal embedding.  

To be sure that this chain does indeed appear in $[u,w]_d$ we look to part b of the definition of a $d$-normal embedding.  The only way a decreasing maximal chain will appear in $[u,w]$, but not $[u,w]_d$ is if it passes through an element that has a run of at least $d+1$ $b$'s.  Since we are deleting elements from right to left, the first time an element with a run of at least $d+1$ elements appears is when an $a$ that is preceded by a $b$ is deleted, and this $a$ is followed by a run of $d$ $b$'s.  This run of $d$ $b$'s must be a run in $u$, so by part b of the definition, $a$ must either be supported, in which case it would not have been deleted, or it must have been the first letter of $w$ in which case it is not preceded by a $b$. 

Conversely, if $\iota =\{i_1,i_2,\dots,i_k\}$ is a $d$-normal embedding of $u$ in $w$ then there is a decreasing maximal chain corresponding to $\iota$ given by deleting the elements $w_j$, $j\in[n]-\iota$, from right to left.   $\square$

The definition of $d$-normal is borne out of the embeddings that give decreasing maximal chains in the CL-shelling of the poset.  Part a of the definition is well understood from the work of Bj\"{o}rner.  We will focus on part b of the definition.  Let $u$ be a word with a run of $d$ $b$'s.  Let $\iota$ be an embedding of $u$ in $w$.  Let $i_j,i_{j+1},\dots,i_{j+d-1}\in\iota$ be the embedding of this run of $d$ $b$'s in $w$.  If we are deleting elements of $w$ from right to left until we are left with $u$ to create a decreasing maximal chain then at some point the remaining letters of $w$ starting with $w_{i_j}$ are all elements of $u$ and we must delete something to the left of $w_{i_j}$ (unless $i_j=1$).  Note that $w_{i_j}$ need not be in a run of $d$ $b$'s in $w$.

There are two possibilities.  Either this run of $d$ $b$'s in $u$ is at the beginning of $u$ or not.  If the run is at the beginning of $u$ then we delete the remaining elements of $w$ one at a time from right to left.  If any of these elements is a $b$ then a string of $d+1$ $b$'s is formed, and hence no decreasing maximal chain corresponds to this embedding.  Thus, there may only be an $a$ preceding this $b$.  

If the run of $d$ $b$'s is not at the beginning of $u$ then in $u$ an $a$ immediately precedes the first $b$ in this run. This $a$ must correspond to $w_{i_j-1}$, otherwise we have a run of unsupported $a$'s or there is a $b$ between this $a$ and $w_{i_j}$ which will eventually create a run of $d+1$ $b$'s.  

\section{Rationality}

We now consider the formal power series algebra $\mathbb{Z}\left\langle\left\langle A\right\rangle\right\rangle$ in the noncommuting variables of $A$ over $\mathbb{Z}$.  Every $f\in\mathbb{Z}\left\langle \left\langle A\right\rangle\right\rangle$ is of the form $$f=\sum_{w\in A^*}c_ww,$$ where $c_w\in\mathbb{Z}$.  We let $f^*=\epsilon+f+f^2+\dots=(\epsilon-f)^{-1}$ for any series $f\in\mathbb{Z}\left\langle \left\langle A\right\rangle\right\rangle$ with no constant term.  Let $f^+=f^*-\epsilon$.  A series $f\in\mathbb{Z}\left\langle \left\langle A\right\rangle\right\rangle$ is called rational if it can be constructed from a finite number of monomials under a finite number of the usual algebraic operations in $\mathbb{Z}\left\langle \left\langle A\right\rangle\right\rangle$ and the $*$ operation.  Clearly, $f^+$ is rational if $f^*$ is.  

We may also consider series of the form $f=\sum_{u\leq w}c_{(u,w)}u\otimes w$ where $u\otimes w$ just represents the ordered pair $(u,w)\in A^*\times A^*$.  Rationality of such a series is defined similarly.  

In~\cite{BjornerReutenauer}, Bj\"{o}rner and Reutenauer showed that the following four series are rational:  
\begin{eqnarray*}
Z(u)&=&\sum_{w\in A^*}\zeta(u,w)w,\\
M(u)&=&\sum_{w\in A^*}\mu(u,w)w,\\
Z_\otimes&=&\sum_{w\in A^*}\zeta(u,w)u\otimes w, \mbox{ and}\\
M_\otimes&=&\sum_{w\in A^*}\mu(u,w)u\otimes w.\\
\end{eqnarray*}

In this section we will use methods similar to those used by Bj\"{o}rner and Sagan~\cite{SagBjorncomposet} to do the same for the series $Z_d(u)$, $M_d(u)$, $Z^d_\otimes$, and $M^d_\otimes$, where these series are defined exactly the same way as those above replacing $A^*$ by $A^*_d$ in each.  In the remainder of this section we will assume that $d=3$.  The reason that we let $d=3$ is to avoid cumbersome formulas and definitions.  Everything that is used to prove these facts for $d=3$ has an obvious generalization.  

We begin with $Z_3(u)$ and $M_3(u)$.  A function $f:S^*\rightarrow T^*$ where $S$ and $T$ are finite sets is {\it multiplicative} if for $u=u_1u_2\dots u_m\in S^*$, we have that $f(u)=f(u_1)f(u_2)\cdots f(u_m)$.  Let $A_3=\{a,ab,abb,abbb\}$.  Let $B=\epsilon+b+bb+bbb$, and notice that $B(A_3)^*=A_3^*$.  Each word $u\in A_3^*$ can be broken uniquely into its maximal runs of $a$'s and $b$'s.  Define the multiplicative function $$z:A_3^*\rightarrow\mathbb{Z}\left\langle \left\langle A_3\right\rangle\right\rangle,$$ by

\begin{eqnarray*}
z(a^k)&=&(aB)^{k-1}a,\\
z(b)&=&(B-\epsilon)a^*,\\
z(bb)&=&(B-\epsilon)a^+ba^*+(B-b-\epsilon)a^*,\mbox{ and}\\
z(bbb)&=&(B-\epsilon)a^+ba^*ba^*+(B-b-\epsilon)a^+ba^*+bbba^*.\\
\end{eqnarray*}

If a run of $a$'s is at the end of the word then let $z(a^k)=(A_3)^k$ for this run of $a$'s.  Let $p_z(u)$ be the {\it prefix} of $Z(u)$ where, $$p_z(u)=\left\{\begin{array}{ll}A_3^*&u\:\mathrm{begins\:with\:}a,\\
(A_3^*a+\epsilon)&u\:\mathrm{begins\:with\:}b.\\ \end{array}\right.$$

Define the multiplicative function $$m:A_3^*\rightarrow\mathbb{Z}\left\langle \left\langle A_3\right\rangle\right\rangle,$$ by

\begin{eqnarray*}
m(a^k)&=&(ab)^*(\epsilon-a)a\left[(\epsilon-b(ab)^*(\epsilon-a))a\right]^{k-1}\\
m(b)&=&\left(\epsilon-b(ab)^*(\epsilon-a)\right)b,\\
m(bb)&=&\left(\epsilon-b(ab)^*(\epsilon-a)\right)b(ab)^*(\epsilon-a)b,\mbox{ and}\\
m(bbb)&=&b\left((ab)^*(\epsilon-a)b\right)^2.\\
\end{eqnarray*}

Let $p_m(u)$ be the {\it prefix} of $M(u)$ where, $$p_m(u)=\left\{\begin{array}{ll}(\epsilon-b)&\mbox{if }u\:\mathrm{begins\:with\:}a,\\
(\epsilon-a)&\mbox{if }u\:\mathrm{begins\:with\:}\:b,bb\:\mathrm{or}\:bbb.\\ \end{array}\right.$$

Let $s_m(u)$ be the {\it suffix} of $M(u)$ where, 
$$s_m(u)=\left\{\begin{array}{ll}(ba)^*(\epsilon-b)&\mbox{if }u\:\mbox{ends in } a,\\
(ab)^*(\epsilon-a)&\mbox{if }u\:\mbox{ends in }b.\end{array}\right.$$

\begin{lem} For any $u\in A^*_3$ we have that $$Z_3(u)=p_z(u)z(u),$$ and $$M_3(u)=p_m(u)m(u)s_m(u).$$ 
\end{lem}

\noindent {\bf Proof:}  We begin by proving the statement for $Z_3(u)$.  We must show that the functions $p_z$ and $z$ will produce each word that contains $u$ exactly once.  This is done by showing that for $w\geq u$, $z(u)$ produces the right-most embedding of $u$ in $w$.  We will explain how this works for $z(a^k)$ and $z(b^\ell)$, how $z(a^k)$ interacts with $z(b^{\ell})$ and how the prefix works.  

The last $a$ in the run $a^k$ is given by the $a$ at the end of $z(a^k)$.  This is clearly under the right-most possible $a$ if we just focus on $z(a^k)$.  The preceding $k-1$ $a$'s are each followed by $B=\epsilon+b+bb+bbb$, so that each word produced by $z(a^k)$ contains $a^k$ exactly once.  Notice that any word produced by $z(a^k)$ begins with a supported $a$ and ends with a supported $a$.   

We now turn our focus to $z(b^\ell)$.  We'll discuss $z(bb)$ as the others are similar.  The first term in $z(bb)$ produces words where the two supported $b$'s are separate, achieved by the $a^+$.  At the beginning of each term in $z(bb)$ we can have a run of $b$'s.  The last $b$ (first term) and the last two $b$'s (second term) are the supported $b$'s.  Finally at the end of each term we have $a^*$ because we can follow our last supported $b$ by as many $a$'s as we like and still maintain a right-most embedding.  

The series $z(a^k)$ and $z(b^\ell)$ interact in the following way.  The last $a$ in $z(a^k)$ is always supported, so if a run of $b$'s follows a run of $a$'s in $u$ then a right-most embedding is still maintained.  The last supported $b$ in $z(b^\ell)$ is followed by $a^*$, so if a run of $a$'s follows a run of $b$'s in $u$ then we maintain a right most embedding.  It is worth mentioning that we have not violated the fact that a run of $b$'s may be of length at most three.

We may essentially place any word from $A_3^*$ at the beginning of any word created by $z$.  We do this with the prefix $p_z(u)$.  If $u$ begins with $a$ we don't have to worry about avoiding a run of three $b$'s.  If $u$ begins with $b$ we do.  It's not hard to see that $p_z(u)$ handles these cases appropriately.  

We will give a similar explanation for $M_3(u)$.  We must show that $p_m(u)m(u)s_m(u)$ will produce the word $w$ exactly $(-1)^{|w|-|u|}{w\choose u}_{3n}$ times.  First we explain $m(u)$.  

The supported $a$'s in the words produced by $m(a^k)$ appear before the left square bracket and the right square bracket.  Before the first supported $a$ we may have a run of alternating $a$'s and $b$'s given by $(ab)^*$.  This run may end in $a$ or $b$ given by $\epsilon-a$.  The $a$ is negative in this term to represent the fact that the parity of $u$ and the word produced $w$ is different if $a$ is placed before the first supported $a$.  A similar argument explains the factor in preceding each of the remaining $k-1$ $a$'s.  Notice that in this case a run of alternating $a$'s and $b$'s must begin with $b$.  This is to maintain appropriately supported repetition sets forcing the embedding to be 3-normal.  

Now, we turn our focus to $m(bbb)$.  The first $b$ and the $b$ preceding the last right parenthesis are the supported $b$'s in any word produced by $m(bbb)$.  Notice that there is nothing preceding the first supported $b$.  This forces the letter preceding the first supported $b$ to be a supported $a$ if $bbb$ is not at the beginning of $u$ satisfying part b of the definition of $d$-normal.  The first supported $b$ can be followed by a run of alternating $a$'s and $b$'s.  The cases $m(b)$ and $m(bb)$ have similar explanations.  

Except in the case when $u$ begins with $bbb$ the prefix and suffix for $m$ merely produce a run of alternating $a$'s and $b$'s at the beginning or end of the produced word taking care of supported repetition sets and preventing any word from being produced more than required.  In the case when $u$ begins with $bbb$, the prefix appends an $a$ or nothing to the beginning in accordance with part b of the definition of $d$-normal.

From our description any word produced by $p_m(u)m(u)s_m(u)$ must give a unique 3-normal embedding of $u$.  This shows that $M_3(u)=p_m(u)m(u)s_m(u)$. $\square$ 

\medskip

The fact that $p_z(u)$, $z(u)$, $p_m(u)$, $m(u)$, and $s_m(u)$ are rational for each $u\in A_3^*$ proves the following theorem.

\begin{thm}  The series $Z_3(u)$ and $M_3(u)$ are rational.  $\square$ \end{thm}

The techniques used above to prove the rationality of $Z_3(u)$ and $M_3(u)$ are a bit cumbersome.  We will use finite state automata to prove that $Z_{\otimes}^3$ and $M_{\otimes}^3$ are rational.  

Let $S$ be an alphabet.  A {\it finite state automaton} is a digraph $D$, with vertex set $V$ and arc set $E$, allowing loops and multi-arcs.  There are unique vertices $\alpha$ and $\omega$, where $\alpha$ is the initial vertex and $\omega$ is the final vertex.  Each arc $e\in E$ is labeled by a monomial $f(e)\in\mathbb{Z}\left\langle \left\langle S\right\rangle\right\rangle$.  A finite walk $W$ with arcs $e_1,e_2,\dots,e_k$ is given the monomial label $$f(W)=\prod_{i=1}^k f(e_i).$$  The series {\it accepted} by $D$ is $$f(D)=\sum_Wf(W),$$ where the sum is over all walks in $D$ from $\alpha$ to $\omega$.  

If $e_1,\dots e_j$ are all arcs from one vertex to another, replacing them by a single arc $e$ and labeling this arc $$\sum_{i=1}^j f(e_i)$$ does not change the series accepted by $D$.  For simplification we will use this procedure.  

It is a well-known fact that a series is rational if and only if it is accepted by a finite state automaton~\cite{BerstelReutenauer}.  We will construct finite state automata accepting $Z_{\otimes}^3$ and $M_{\otimes}^3$ to prove the following theorem for $d=3$.  The pattern in the automata will be obvious and generalizable.  

\vspace{.25in}

\begin{picture}(400,450)

\put(200,400){\circle{20}}
\put(120,400){\circle{20}}
\put(280,400){\circle{20}}
\put(195,400){$\alpha_2$}
\put(115,400){$\alpha_3$}
\put(275,400){$\alpha_1$}
\put(270,400){\vector(-1,0){60}}
\put(190,400){\vector(-1,0){60}}
\put(130,395){\vector(2,-3){57}}
\put(270,395){\vector(-2,-3){57}}
\put(200,390){\vector(0,-1){78}}
\put(213,300){\vector(2,3){59}}
\put(150,405){$\epsilon\otimes b$}
\put(230,405){$\epsilon\otimes b$}
\put(190,340){\rotatebox{90}{$\epsilon\otimes a$}}
\put(134,365){\rotatebox{305}{$\epsilon\otimes a$}}
\put(224,345){\rotatebox{55}{$\epsilon\otimes a$}}
\put(244,335){\rotatebox{55}{$\epsilon\otimes b$}}

\put(200,300){\circle{25}}
\put(195,297){{\large $\alpha$}}
\put(190,145){\circle{20}}
\put(185,142){$\alpha_4$}
\put(50,180){\circle{20}}
\put(350,180){\circle{20}}
\put(100,50){\circle{20}}
\put(300,50){\circle{20}}
\put(275,140){\circle{20}}
\put(45,175){$\beta_1$}
\put(345,175){$\gamma_2$}
\put(270,135){$\gamma_3$}
\put(95,45){$\beta_2$}
\put(295,45){$\beta_3$}
\put(300,230){\circle{20}}
\put(295,225){$\gamma_1$}
\put(190,290){\vector(-4,-3){135}}
\put(50,170){\vector(1,-2){55}}
\put(110,50){\vector(1,0){180}}
\put(345,170){\vector(-1,-3){37}}
\put(305,220){\vector(1,-1){33}}
\put(212,298){\vector(4,-3){80}}
\put(60,180){\vector(4,-1){115}}
\put(175,144){\vector(-4,1){115}}
\put(195,290){\vector(0,-1){135}}
\put(105,60){\vector(1,1){75}}
\put(290,50){\vector(-1,1){90}}
\put(265,140){\vector(-1,0){60}}
\put(340,180){\vector(-2,-1){60}}
\put(200,145){\vector(1,1){85}}
\put(288,225){\vector(-1,-1){85}}
\put(340,183){\line(-2,1){62}}
\put(245,173){\rotatebox{45}{$a\otimes a$}}

\put(188,295){\line(-1,0){30}}
\put(158,305){\vector(1,0){30}}
\put(158,300){\oval(10,10)[l]}

\put(310,30){\oval(35,30)[bl]}
\put(310,30){\oval(35,30)[r]}
\put(292.5,30){\vector(0,1){10}}

\put(90,30){\oval(35,30)[br]}
\put(90,30){\oval(35,30)[l]}
\put(107.5,30){\vector(0,1){10}}

\put(35,165){\oval(30,30)[br]}
\put(35,165){\oval(30,30)[l]}
\put(50,165){\vector(0,1){5}}

\put(185,136){\line(0,-1){30}}
\put(195,106){\vector(0,1){30}}
\put(190,106){\oval(10,10)[b]}

\put(310,230){\line(1,0){80}}
\put(390,230){\line(0,-1){230}}
\put(390,0){\line(-1,0){240}}
\put(150,0){\vector(-1,1){40}}

\put(126,298){$\epsilon\otimes a$}

\put(315,220){\rotatebox{310}{$\epsilon\otimes b$}}
\put(200,5){$b\otimes b$}
\put(300,5){$\epsilon\otimes a$}
\put(250,280){\rotatebox{320}{$\epsilon\otimes b$}}
\put(200,230){\rotatebox{270}{$a\otimes a$}}
\put(230,185){\rotatebox{55}{$\epsilon\otimes b$}}
\put(220,145){$a\otimes a$}
\put(295,162){\rotatebox{28}{$\epsilon\otimes b$}}
\put(315,110){\rotatebox{78}{$b\otimes b$}}
\put(177,93){$a\otimes a$}
\put(200,55){$b\otimes b$}
\put(100,230){\rotatebox{42}{$b\otimes b$}}
\put(75,5){$\epsilon\otimes a$}
\put(110,150){\rotatebox{345}{$b\otimes b$}}
\put(110,170){\rotatebox{345}{$a\otimes a$}}
\put(20,140){$\epsilon\otimes a$}
\put(238,90){\rotatebox{315}{$a\otimes a$}}
\put(78,120){\rotatebox{295}{$b\otimes b$}}
\put(130,100){\rotatebox{45}{$a\otimes a$}}
\put(0,450){{\bf Figure 2: $Z_{\otimes}^3$}}

\end{picture}

\vspace{.5in}
  
\begin{picture} (500,700)

\put(0,700){{\bf Figure 3: $M_{\otimes}^3$}}
\put(100,490){\circle{20}}
\put(200,490){\circle{20}}
\put(195,490){$\alpha_3$}
\put(95,490){$\alpha_4$}
\put(50,400){\circle{20}}
\put(160,400){\circle{20}}
\put(335,400){\circle{20}}
\put(45,400){$\beta_1$}
\put(155,400){$\gamma_1$}
\put(330,400){$\delta_1$}
\put(50,300){\circle{20}}
\put(50,200){\circle{20}}
\put(50,100){\circle{20}}
\put(45,300){$\beta_2$}
\put(45,200){$\beta_3$}
\put(45,100){$\beta_4$}
\put(130,300){\circle{20}}
\put(190,300){\circle{20}}
\put(130,200){\circle{20}}
\put(190,200){\circle{20}}
\put(130,100){\circle{20}}
\put(190,100){\circle{20}}
\put(125,300){$\gamma_2$}
\put(125,200){$\gamma_3$}
\put(125,100){$\gamma_4$}
\put(185,300){$\gamma_5$}
\put(185,200){$\gamma_6$}
\put(185,100){$\gamma_7$}
\put(260,300){\circle{20}}
\put(260,200){\circle{20}}
\put(260,100){\circle{20}}
\put(335,300){\circle{20}}
\put(335,200){\circle{20}}
\put(335,100){\circle{20}}
\put(405,300){\circle{20}}
\put(405,200){\circle{20}}
\put(405,100){\circle{20}}
\put(255,300){$\delta_2$}
\put(255,200){$\delta_3$}
\put(255,100){$\delta_4$}
\put(330,300){$\delta_5$}
\put(330,200){$\delta_6$}
\put(330,100){$\delta_7$}
\put(400,300){$\delta_8$}
\put(400,200){$\delta_9$}
\put(400,100){$\delta_{10}$}
\put(210,510){\oval(35,30)[tl]}
\put(210,510){\oval(35,30)[r]}
\put(192.5,510){\vector(0,-1){10}}
\put(190,495){\vector(-1,0){80}}
\put(110,485){\vector(1,0){78}}
\put(140,500){$-\epsilon\otimes b$}
\put(140,475){$-\epsilon\otimes a$}
\put(90,490){\vector(-1,-2){40}}
\put(200,480){\vector(-1,-2){35}}
\put(110,480){\vector(2,-3){45}}
\put(202,480){\vector(3,-2){120}}
\put(50,390){\vector(0,-1){80}}
\put(50,290){\vector(0,-1){80}}
\put(45,190){\vector(0,-1){80}}
\put(55,110){\vector(0,1){80}}
\put(155,390){\vector(-1,-3){25}}
\put(165,390){\vector(1,-3){25}}
\put(130,290){\vector(0,-1){80}}
\put(125,190){\vector(0,-1){80}}
\put(135,110){\vector(0,1){80}}
\put(190,290){\vector(0,-1){80}}
\put(185,190){\vector(0,-1){80}}
\put(195,110){\vector(0,1){80}}
\put(330,390){\vector(-2,-3){59}}
\put(335,390){\vector(0,-1){80}}
\put(340,390){\vector(2,-3){57}}
\put(260,290){\vector(0,-1){80}}
\put(255,190){\vector(0,-1){80}}
\put(265,110){\vector(0,1){80}}
\put(335,290){\vector(0,-1){80}}
\put(330,190){\vector(0,-1){80}}
\put(340,110){\vector(0,1){80}}
\put(405,290){\vector(0,-1){80}}
\put(400,190){\vector(0,-1){80}}
\put(410,110){\vector(0,1){80}}
\put(160,248){\vector(2,3){27}}
\put(140,100){\line(1,0){20}}
\put(160,100){\line(0,1){148}}
\put(140,200){\line(1,0){20}}
\put(270,200){\line(1,0){20}}
\put(270,100){\line(1,0){20}}
\put(290,100){\line(0,1){150}}
\put(290,250){\vector(1,1){40}}
\put(345,200){\line(1,0){20}}
\put(345,100){\line(1,0){20}}
\put(365,100){\line(0,1){150}}
\put(365,250){\vector(1,1){40}}
\put(250,200){\line(-1,0){20}}
\put(250,100){\line(-1,0){20}}
\put(230,200){\line(0,-1){120}}
\put(230,80){\line(1,0){205}}
\put(435,80){\line(0,1){170}}
\put(435,250){\vector(-3,4){30}}
\put(200,530){$a\otimes a$}
\put(40,335){\rotatebox{90}{$b\otimes b$}}
\put(40,235){\rotatebox{90}{$-\epsilon\otimes a$}}
\put(35,135){\rotatebox{90}{$-\epsilon\otimes b$}}
\put(55,135){\rotatebox{90}{$-\epsilon\otimes a$}}
\put(48,430){\rotatebox{65}{$a\otimes a+\epsilon\otimes\epsilon$}}
\put(125,345){\rotatebox{65}{$b\otimes b$}}
\put(180,360){\rotatebox{295}{$bb\otimes bb$}}
\put(120,235){\rotatebox{90}{$-\epsilon\otimes a$}}
\put(115,135){\rotatebox{90}{$-\epsilon\otimes b$}}
\put(135,135){\rotatebox{90}{$-\epsilon\otimes a$}}
\put(160,185){\rotatebox{90}{$b\otimes b$}}
\put(192,235){\rotatebox{90}{$-\epsilon\otimes a$}}
\put(175,135){\rotatebox{90}{$-\epsilon\otimes b$}}
\put(200,135){\rotatebox{90}{$-\epsilon\otimes a$}}
\put(105,465){\rotatebox{305}{$a\otimes a+\epsilon\otimes\epsilon$}}
\put(182,435){\rotatebox{65}{$\epsilon\otimes \epsilon$}}
\put(250,450){\rotatebox{330}{$a\otimes (a-ba)$}}
\put(275,330){\rotatebox{60}{$b\otimes b$}}
\put(325,330){\rotatebox{90}{$bb\otimes bb$}}
\put(365,360){\rotatebox{305}{$bbb\otimes bbb$}}
\put(250,235){\rotatebox{90}{$-\epsilon\otimes a$}}
\put(245,135){\rotatebox{90}{$-\epsilon\otimes b$}}
\put(268,135){\rotatebox{90}{$-\epsilon\otimes a$}}
\put(320,70){$bb\otimes bb$}
\put(290,190){\rotatebox{90}{$b\otimes b$}}
\put(325,235){\rotatebox{90}{$-\epsilon\otimes a$}}
\put(320,135){\rotatebox{90}{$-\epsilon\otimes b$}}
\put(340,135){\rotatebox{90}{$-\epsilon\otimes a$}}
\put(395,235){\rotatebox{90}{$-\epsilon\otimes a$}}
\put(390,135){\rotatebox{90}{$-\epsilon\otimes b$}}
\put(410,135){\rotatebox{90}{$-\epsilon\otimes a$}}
\put(370,190){\rotatebox{90}{$b\otimes b$}}
\put(190,495){\line(-1,1){20}}
\put(170,515){\line(-1,0){125}}
\put(45,515){\vector(0,-1){105}}
\put(80,520){$\epsilon\otimes\epsilon$}

\put(315,610){\circle{20}}
\put(310,610){$\alpha_3$}
\put(345,610){\circle{20}}
\put(340,610){$\beta_2$}
\put(385,610){\circle{20}}
\put(380,610){$\gamma_2$}
\put(415,610){\circle{20}}
\put(410,610){$\gamma_5$}
\put(366,685){\circle{25}}
\put(360,683){{\large $\alpha_2$}}
\put(415,685){\line(-1,0){37}}
\put(315,685){\line(1,0){38}}
\put(315,685){\vector(0,-1){65}}
\put(345,685){\vector(0,-1){65}}
\put(385,685){\vector(0,-1){65}}
\put(415,685){\vector(0,-1){65}}
\put(302,630){\rotatebox{90}{$a\otimes a$}}
\put(332,630){\rotatebox{90}{$b\otimes b$}}
\put(372,630){\rotatebox{90}{$b\otimes b$}}
\put(402,630){\rotatebox{90}{$bb\otimes bb$}}

\put(125,610){\circle{20}}
\put(120,610){$\alpha_3$}
\put(155,610){\circle{20}}
\put(150,610){$\beta_2$}
\put(195,610){\circle{20}}
\put(190,610){$\gamma_2$}
\put(225,610){\circle{20}}
\put(220,610){$\gamma_5$}
\put(175,685){\circle{25}}
\put(170,683){{\large $\alpha_1$}}
\put(125,685){\line(1,0){39}}
\put(225,685){\line(-1,0){39}}
\put(125,685){\vector(0,-1){65}}
\put(155,685){\vector(0,-1){65}}
\put(195,685){\vector(0,-1){65}}
\put(225,685){\vector(0,-1){65}}
\put(112,630){\rotatebox{90}{$a\otimes a$}}
\put(142,630){\rotatebox{90}{$b\otimes b$}}
\put(182,630){\rotatebox{90}{$b\otimes b$}}
\put(212,630){\rotatebox{90}{$bb\otimes bb$}}

\put(15,675){\circle{25}}
\put(10,675){{\large$\alpha$}}
\put(15,662){\vector(0,-1){72}}
\put(15,580){\circle{20}}
\put(10,580){$\alpha_1$}
\put(15,515){\circle{20}}
\put(10,515){$\alpha_2$}
\put(12,570){\vector(0,-1){45}}
\put(18,525){\vector(0,1){45}}
\put(2,600){\rotatebox{90}{$\epsilon\otimes (\epsilon-a)$}}
\put(0,530){\rotatebox{90}{$-\epsilon\otimes b$ }}
\put(19,530){\rotatebox{90}{$-\epsilon\otimes a$ }}
\put(45,540){\circle{20}}
\put(40,535){$\delta_1$}
\put(45,675){\vector(0,-1){125}}
\put(32,570){\rotatebox{90}{$\epsilon\otimes(\epsilon-a)+a\otimes a$}}
\put(45,675){\line(-1,0){19}}

\end{picture}

\begin{thm}  For any $d$, $Z_{\otimes}^d$ and $M_{\otimes}^d$ are rational.\end{thm}    

\noindent {\bf Proof:}  The automata in Figures 2 and 3 are for $Z_{\otimes}^3$ and $M_{\otimes}^3$ respectively.  We will use them to explain why these automata work and how they are generalizable to any $d$.  For clarity we left some arcs off of the diagrams.  In Figure 2, there is an arc labeled $\epsilon\otimes\epsilon$ from each node to $\omega$.  In Figure 3, there is an arc labeled $\epsilon\otimes\epsilon$ from each of $\alpha_1$, $\alpha_2$, $\alpha_3$, $\beta_2$, $\beta_3$, $\beta_4$, $\gamma_5$, $\gamma_6$, $\gamma_7$, $\delta_8$, $\delta_9$, and $\delta_{10}$ to $\omega$.  There is also an arc labeled $a\otimes a$ from each of $\beta_2$, $\beta_3$, $\beta_4$, $\gamma_5$, $\gamma_6$, $\gamma_7$, $\delta_8$, $\delta_9$, and $\delta_{10}$ to each of $\alpha_3$, $\beta_1$, $\gamma_1$, and $\delta_1$.  Also in Figure 3, we separated the nodes $\alpha$, $\alpha_1$ and $\alpha_2$ from the rest of the diagram for clarity. 

We claim that for any $u=u_1u_2\dots u_k\leq w=w_1w_2\dots w_n$ there is a unique walk from $\alpha$ to $\omega$ in the automaton for $Z_{\otimes}^3$ labeled $u\otimes w$.  Consider the right-most embedding $\iota=\{i_1,i_2,\dots,i_k\}$ of $u$ in $w$.  We'll build $u\otimes w$ from this embedding.  The first $i_1-1$ letters in $w$ are not supported, so this portion of $w$ is constructed using the nodes $\alpha$, $\alpha_1$, $\alpha_2$, and $\alpha_3$.  Notice that any word in $A_3^*$ can be built uniquely by using the nodes $\alpha$, $\alpha_1$, $\alpha_2$, and $\alpha_3$.  Now, if $u_1=a$ and $w_{i_1}$ is preceded by an $a$ then we must go along the arc from $\alpha$ to $\alpha_4$.  If $w_{i_1}$ is preceded by a run of $b$'s then we must go along the arc from $\alpha$ to $\gamma_1$, construct this run of $b$'s and then go to $\alpha_4$.  If $u_1=b$ and $w_{i_1}$ is preceded by an $a$ then at this point the walk goes along the arc from $\alpha$ to $\beta_1$.  If $u_1=b$ and $w_{i_1}$ is preceded by a run of $b$'s then the walk must go along the arc from $\alpha$ to $\gamma_1$, build up the preceding run of unsupported $b$'s and then go on to either $\beta_2$ or $\beta_3$ from $\gamma_1$ or $\gamma_2$. 

Now, there is no way to get back to the $\alpha$, $\alpha_1$, $\alpha_2$, and $\alpha_3$ part of the digraph.  Notice that if $u_1=a$ then to construct the unsupported part of $w$ between $w_{i_1}$ and $w_{i_2}$ we can only construct a run of $b$'s given by the cycle $\gamma_1\rightarrow\gamma_2\rightarrow\gamma_3\rightarrow\alpha_4$.  If $u_1=b$ then the unsupported part of $w$ between $w_{i_1}$ and $w_{i_2}$ can only be a run of $a$'s given by the loops labeled $\epsilon\otimes a$ at each of $\beta_1$, $\beta_2$ and $\beta_3$.  This is forcing the embedding $\iota$ to be the right-most embedding.  The construction continues as above for the remainder of $u\otimes w$.  

It is important to note here that the complication of the diagram develops from the fact that we must be careful to avoid producing a run of more than three $b$'s in $u$ or $w$.  The portion of the digraph involving $\beta_1$, $\beta_2$ and $\beta_3$ controls the supported $b$'s in $w$, and the portion involving $\gamma_1$, $\gamma_2$ and $\gamma_3$ controls the unsupported $b$'s.  

Our comments above about the automaton forcing $\iota$ to be the right-most embedding proves that each $u\otimes w$ is produced by a unique walk, and hence the automaton accepts $Z_{\otimes}^3$.  To generalize this to any $Z_{\otimes}^d$, we would merely extend the $\alpha$, $\beta$, and $\gamma$ portions of the digraph appropriately.

We turn our focus to the automaton for $M_\otimes^3$.  We claim that there are ${w\choose u}_{3n}$ walks each of which produces $(-1)^{|w|+|u|}u\otimes w$.  Each word begins with a run of alternating unsupported $a$'s and $b$'s.  This is represented by the $\alpha$, $\alpha_1$, $\alpha_2$ portion of the automaton.  If the embedded word begins with a run of three $b$'s then the automaton goes directly from $\alpha$ to the $\delta$ portion of the automaton.  Once this run has been constructed the walk must move on to its first supported letter.  If the first supported letter is an $a$ then the walk goes from $\alpha_1$ or $\alpha_2$ to $\alpha_3$.  If the first supported letter is a $b$ then the walk goes into the $\beta$, $\gamma$ or $\delta$ portions of the automaton.  

Once the walk leaves the $\alpha$, $\alpha_1$, $\alpha_2$ portion of the automaton it cannot return.  The $\alpha_3$ and $\alpha_4$ portion of the automaton handles runs of supported $a$'s and runs of alternating unsupported $a$'s and $b$'s.

Now, the three legs of the diagram labeled with $\beta$'s, $\gamma$'s and $\delta$'s respectively, are controlling the three possibilities for a run of $b$'s in $u$.  If $u$ has a run of three $b$'s, the walk producing $u\otimes w$ must pass through the $\delta$ part of the diagram.  By part b in the definition of $d$-normal, the first supported $b$ in a run of three $b$'s must either be at the beginning of $w$, handled by the arc from $\alpha$ to $\delta_1$, or the first supported $b$ must be preceded by a supported $a$, handled by the arc from $\alpha_3$ to $\delta_1$.  Each of the three different legs of the $\delta$ portion of the diagram represents a supported $b$.  So once the walk is on any of the nodes $\delta_8$, $\delta_9$, or $\delta_{10}$ all three $b$'s in $u$ have been produced.  Notice that in this section between any two supported $b$'s there can be a run of alternating $a$'s and $b$'s, given by the exchange between the bottom two nodes of each leg.  Also, each arc going into $\delta_1$ is labeled with a supported $a$ immediately preceding the portion of the word to be produce.  This again is assuring that the embedding of any walk is 3-normal. 

The portion marked with $\beta$'s controls situations where $u$ has a run of just one $b$ and the portion marked with $\gamma$'s takes care of two $b$'s.  This shows that every walk from $\alpha$ to $\omega$ in this automaton produces $u\otimes w$ according to a 3-normal embedding.  Thus, each 3-normal embedding $\iota$ of $u$ in $w$ is given by a unique walk from $\alpha$ to $\omega$.  

Finally, notice that anytime a letter in a label on an arc is unsupported the sign of the monomial label changes.  This takes care of the sign of $\mu(u,w)$.  Thus, there are exactly ${w\choose u}_{3n}$ walks producing $(-1)^{|w|+|u|}u\otimes w$.  $\square$  

\section{Generating Functions in Commuting Variables}

We now turn to generating functions in commuting variables.  We will consider generating functions in terms of the norm of the composition $\alpha$, $|\alpha|$.  If $n_k$ is the number of $k$'s in $\alpha$ then we may redefine $|\alpha|=\sum_{k\geq1}n_kk$.  Let $\alpha\leftrightarrow u$, then the {\it type} of $\alpha$ is $\tau(\alpha)=(n_1,n_2,\dots,n_n,r)$, where $r$ is the number of runs of $a$'s in $u$.  

Each time $k$ appears in $u$ replace $k$ by $x^k$, where $x$ is a commuting variable.  Then we obtain the norm generating functions $$Z_d(\alpha;x)=\sum_{\alpha\leq \beta}x^{|\beta|}$$ and $$M_d(\alpha;x)=\sum_{\alpha\leq\beta}\mu(\alpha,\beta)x^{|\beta|}.$$

To avoid cumbersome formulas and because these generalize simply to any $d$ we will focus on the case when $d=3$.  If $\alpha\in C_4$ corresponds to $u\in A_3^*$ then any $a$ that is not immediately followed by a $b$ represents a 1 from $\alpha$ and $b$, $bb$, and $bbb$ represent 2,3, and 4 respectively.  Let $[k]_x$ be the polynomial $1+x+\dots+x^{k-1}$.

Let $p_z(u;x)$ and $z(u;x)$ be the formal power series in $\mathbb{Z}\left\langle \left\langle x\right\rangle\right\rangle$ obtained from $z(u)$ and $p_z(u)$ respectively by replacing each letter in $p_z(u)$ and $z(u)$ by $x$ and multiplying the whole word by $x$.  We have that $Z_3(u;x)=xp_z(u;x)z(u;x)$.  Defining $p_m(u;x)$, $m(u;x)$ and $s_m(u;x)$ in a similar way gives us that $M_3(u;x)=xp_m(u;x)m(u;x)s_m(u;x)$.  Note the $x$ at the beginning of these takes care of the fact that we dropped the initial $a$ from $\phi(\alpha)$.  

Suppose $u$ has type $\tau(u)=(n_1,n_2,n_3,n_4,r)$ then Lemma 4.1 gives us that $Z_3(u;x)$ is one of the following rational functions.  The first two correspond to $u$ beginning with $a$ and the last two correspond to $u$ beginning with $b$.  The first and third correspond to $u$ ending with $a$ and the second and fourth correspond to $u$ ending with $b$.
 
\begin{enumerate}
 \item $x\left(\frac{[4]_x}{1-x[4]_x}\right)\left(x[4]_x\right)^{n_1-r+1}\left(\frac{x[3]_x}{1-x}\right)^{n_2}\left(\frac{x^3[3]_x}{(1-x)^2}+\frac{x^2[2]_x}{1-x}\right)^{n_3}\left(\frac{x^4[3]_x}{(1-x)^3}+\frac{x^4[2]_x}{(1-x)^2}+\frac{x^3}{1-x}\right)^{n_4}$

\item $x\left(\frac{[4]_x}{1-x[4]_x}\right)\left(x[4]_x\right)^{n_1-r}\left(\frac{x[3]_x}{1-x}\right)^{n_2}\left(\frac{x^3[3]_x}{(1-x)^2}+\frac{x^2[2]_x}{1-x}\right)^{n_3}\left(\frac{x^4[3]_x}{(1-x)^3}+\frac{x^4[2]_x}{(1-x)^2}+\frac{x^3}{1-x}\right)^{n_4}$

\item $x\left(\frac{x[4]_x}{1-x[4]_x}+1\right)\left(x[4]_x\right)^{n_1-r+1}\left(\frac{x[3]_x}{1-x}\right)^{n_2}\left(\frac{x^3[3]_x}{(1-x)^2}+\frac{x^2[2]_x}{1-x}\right)^{n_3}\left(\frac{x^4[3]_x}{(1-x)^3}+\frac{x^4[2]_x}{(1-x)^2}+\frac{x^3}{1-x}\right)^{n_4}$

\item $x\left(\frac{x[4]_x}{1-x[4]_x}+1\right)\left(x[4]_x\right)^{n_1-r}\left(\frac{x[3]_x}{1-x}\right)^{n_2}\left(\frac{x^3[3]_x}{(1-x)^2}+\frac{x^2[2]_x}{1-x}\right)^{n_3}\left(\frac{x^4[3]_x}{(1-x)^3}+\frac{x^4[2]_x}{(1-x)^2}+\frac{x^3}{1-x}\right)^{n_4}$

\end{enumerate}
Again, if $u$ has type $\tau(u)=(n_1,n_2,n_3,n_4,r)$ then $M_3(u;x)$ is one of the following.  Formula 5 corresponds to $u$ ending with $a$ and formula 6 corresponds to $u$ ending with $b$.    

\begin{enumerate}
\setcounter{enumi}{4}
\item $x(1-x)\left(\frac{1}{(1+x)^r}\right)\left(1-\frac{x}{1+x}\right)^{n_1-r}x^{n_1-r+1}\left(1-\frac{x}{1+x}\right)^{n_2}\left(\frac{x^2}{1+x}(1-\frac{x}{1+x})\right)^{n_3}\left(\frac{x^3}{(1+x)^2}\right)^{n_4}\frac{1}{1+x}$

\item $x(1-x)\left(\frac{1}{(1+x)^r}\right)\left(1-\frac{x}{1+x}\right)^{n_1-r}x^{n_1-r}\left(1-\frac{x}{1+x}\right)^{n_2}\left(\frac{x^2}{1+x}(1-\frac{x}{1+x})\right)^{n_3}\left(\frac{x^3}{(1+x)^2}\right)^{n_4}\frac{1}{1+x}$

\end{enumerate}

The following theorem is now an immediate consequence of Lemma 4.1.  

\begin{thm}  The norm generating functions $Z_3(u;x)$ and $M_3(u;x)$ for $u$ with type $\tau(u)=(n_1,n_2,n_3,n_4,r)$ are as stated above.  $\square$ \end{thm}

The author did not get a chance to explore whether there were nice generating functions for powers of $\zeta$.  These would be interesting because $\zeta^m(u,w)$ is the number of chains of length $m$ beginning with $u$ and ending with $w$. 

We would to thank Anders Bj\"{o}rner and Richard Stanley without whose prior work this paper would not have been possible.  We would also like to thank the referees for making this paper much easier to read.

\bibliography{Bib}
\end{document}